\documentclass[conference]{IEEEtran}
\IEEEoverridecommandlockouts

\usepackage{cite}
\usepackage{amsmath,amssymb,amsfonts}
\usepackage{algorithmic}
\usepackage{graphicx}
\usepackage{textcomp}
\usepackage{xcolor}
\def\BibTeX{{\rm B\kern-.05em{\sc i\kern-.025em b}\kern-.08em
    T\kern-.1667em\lower.7ex\hbox{E}\kern-.125emX}}

\usepackage{glossaries}

\newcommand{\eeq}{\end{equation}}

\newcommand{\br}{\mbox{\boldmath $r$}}
\newcommand{\bT}{\mbox{\boldmath $T$}}

\newcommand{\bp}{\mbox{\boldmath $p$}}

\newcommand{\bzero}{\mbox{\boldmath $0$}}
\newcommand{\bM}{\mbox{\boldmath $M$}}

\newcommand{\bA}{\mbox{\boldmath $A$}}

\newcommand{\bC}{\mbox{\boldmath $C$}}
\newcommand{\bX}{\mbox{\boldmath $X$}}
\newcommand{\bh}{\mbox{\boldmath $h$}}
\newcommand{\bH}{\mbox{\boldmath $H$}}

\newcommand{\bc}{\mbox{\boldmath $c$}}

\newcommand{\bgamma}{\mbox{\boldmath $\gamma$}}

\newcommand{\bU}{\mbox{\boldmath $U$}}
\newcommand{\bGamma}{\mbox{\boldmath $\Gamma$}}

\newcommand{\bG}{\mbox{\boldmath $G$}}

\newcommand{\bW}{\mbox{\boldmath $W$}}

\newcommand{\bu}{\mbox{\boldmath $u$}}

\newcommand{\bn}{\mbox{\boldmath $n$}}

\newcommand{\bR}{\mbox{\boldmath $R$}}

\newcommand{\bI}{\mbox{\boldmath $I$}}

\newcommand{\by}{\mbox{\boldmath $y$}}
\newcommand{\ds}{\displaystyle}
\newcommand{\bw}{\mbox{\boldmath $w$}}

\newcommand{\beq}{\begin{equation}}

\newcommand{\tr}      {{\mathrm{tr}}}    

\makeglossaries

\newacronym{mac}{MAC}{multiple-access channel}
\newacronym{bc}{BC}{broadcast channel}
\newacronym{mimo}{MIMO}{multiple-input multiple-output}
\newacronym{siso}{SISO}{single-input single-output}
\newacronym{sc}{SC}{single-carrier}
\newacronym{mc}{MC}{multi-carrier}
\newacronym{ofdma}{OFDMA}{orthogonal frequency division multiple access}
\newacronym{af}{AF}{amplify-and-forward}
\newacronym{df}{DF}{decode-and-forward}
\newacronym{cf}{CF}{compress-and-forward}
\newacronym{mwrc}{MWRC}{multi-way relay channel}
\newacronym{pde}{PDE}{partial data exchange}
\newacronym{fde}{FDE}{full data exchange}
\newacronym{iid}{i.i.d.\@}{independent and identically distributed}
\newacronym{awgn}{AWGN}{additive white Gaussian noise}
\newacronym{awg}{AWG}{additive white Gaussian}
\newacronym{sic}{SIC}{successive interference cancellation}
\newacronym{snr}{SNR}{signal-to-noise ratio}
\newacronym{sinr}{SINR}{signal-to-interference-plus-noise ratio}
\newacronym{ber}{BER}{bit error rate}
\newacronym{zf}{ZF}{zero-forcing}
\newacronym{mmse}{MMSE}{minimum mean square error}
\newacronym{sud}{SUD}{single user decoding}
\newacronym{dof}{DoF}{degrees of freedom}
\newacronym{gdof}{GDoF}{generalized degrees of freedom}
\newacronym{nnc}{NNC}{noisy network coding}
\newacronym{dmn}{DMN}{discrete memoryless network}
\newacronym{csi}{CSI}{channel state information}
\newacronym{ee}{EE}{energy efficiency}
\newacronym{ian}{IAN}{treating interference as noise}
\newacronym{snd}{SND}{simultaneous non-unique decoding}
\newacronym{brd}{BRD}{best response dynamics}
\newacronym{br}{BR}{best response}
\newacronym{ne}{NE}{Nash equilibrium}
\newacronym{lhs}{LHS}{left-hand side}
\newacronym{rhs}{RHS}{right-hand side}
\newacronym{gee}{GEE}{global energy efficiency}
\newacronym{wsee}{WSEE}{weighted sum energy efficiency}
\newacronym{wpee}{WPEE}{weighted product energy efficiency}
\newacronym{wmee}{WMEE}{weighted minimum energy efficiency}
\newacronym{kkt}{KKT}{Karush-Kuhn-Tucker}
\newacronym{pc}{PC}{pseudo-concave}
\newacronym{qc}{QC}{quasi-concave}
\newacronym{ql}{QL}{quasi-linear}
\newacronym{pl}{PL}{pseudo-linear}
\newacronym{spc}{SPC}{strictly pseudo-concave}
\newacronym{sqc}{SQC}{strictly quasi-concave}
\newacronym{lfp}{LFP}{linear fractional problem}
\newacronym{clfp}{CLFP}{concave-linear fractional problem}
\newacronym{ccfp}{CCFP}{concave-convex fractional problem}
\newacronym{mmfp}{MMFP}{max-min fractional problem}
\newacronym{sorp}{SoRP}{sum-of-ratios problem}
\newacronym{porp}{PoRP}{product-of-ratios problem}
\newacronym{srp}{SRP}{single-ratio problem}
\newacronym{brb}{BRB}{branch-reduce-and-bound}
\newacronym{qos}{QoS}{quality-of-service}
\newacronym{comp}{CoMP}{cooperative multi-point}
\newacronym{ue}{UE}{user equipment}
\newacronym{bs}{BS}{base station}
\newacronym{mrc}{MRC}{maximum ratio combining}
\newacronym{d2d}{D2D}{device-to-device}
\newacronym{lmmse}{LMMSE}{linear minimum mean square error}
\newacronym{ris}{RIS}{reconfigurable intelligent surface}
\newacronym{svd}{SVD}{singular values decomposition}

\begin{document}

\title{Energy Efficiency in RIS-Aided Wireless Networks: Active or Passive RIS?
\thanks{The work of R. K. Fotock and A. Zappone was supported by the European Commission through the H2020-MSCA-ITN-METAWIRELESS project, grant agreement 956256. The work of M. Di Renzo was supported in part by the European Commission through the H2020 ARIADNE project, grant agreement 871464 and the H2020 RISE-6G project, grant agreement 101017011.}
}

\author{\IEEEauthorblockN{1\textsuperscript{st} Robert K. Fotock}
\IEEEauthorblockA{\textit{CNIT} \\
Cassino, Italy \\
rfotock@cnit.it}
\and
\IEEEauthorblockN{2\textsuperscript{nd} Alessio Zappone}
\IEEEauthorblockA{\textit{CNIT \& Univ. of Cassino and Southern Lazio} \\
Cassino, Italy \\
alessio.zappone@unicas.it}
\and
\IEEEauthorblockN{3\textsuperscript{rd} Marco Di Renzo}
\IEEEauthorblockA{\textit{CentraleSupelec-University} \\
Paris, France \\
marco.di-renzo@universite-paris-saclay.fr}
}

\maketitle

\begin{abstract}
This work addresses the comparison between active and passive RISs in wireless networks, with reference to the system energy efficiency (EE). Two provably convergent and computationally-friendly EE maximization algorithms are developed, which optimize the reflection coefficients of the RIS, the transmit powers, and the linear receive filters. Numerical results show the performance of the proposed methods and discuss the operating points in which active or passive RISs should be preferred from an energy-efficient perspective. 
\end{abstract}

\section{Introduction}
Recently, reconfigurable intelligent surfaces (RISs) have emerged as one of the main technologies for future wireless networks, for their ability to provide many degrees of freedom with limited power consumption \cite{Ref10,RuiZhang_COMMAG,SmartWireless,HuangMag2020}. In particular, energy efficiency (EE) is a major requirement of future wireless networks, also because, according to recent studies, 5G has not achieved the 2000x EE increase that was sought after \cite{David5G}.   

From an energy perspective, the nearly passive behavior of RISs has been recognized as a major advantage, but, on the other hand, it also limits the rate performance that a RIS-aided network can ensure \cite{Ref23,Ref11}, especially when no direct path exists \cite{Ref12,Ref22}. For this reason, recently, the use of active RIS has started to be investigated, i.e. the RIS is equipped with analog amplifiers that allow it to increase the amplitude of the incoming signal \cite{Ref12}. 
The idea of active RIS was influenced by the related field of backscatter communications, where the use of active loads was proposed in \cite{Ref25}. In \cite{Ref24}, the use of reflection-type amplifiers is proposed to implement an active RIS. In this context, a joint transmit beamforming and reflect precoding algorithm is proposed. A comparison between active and passive RIS in terms of rate performance is performed in \cite{Ref16}, considering the optimization of the position of the RIS in both uplink and downlink communications. In \cite{Ref21} a comparative study on channel estimation and spectral efficiency for both purely passive and hybrid RIS architectures was carried out. A hardware architecture for active RIS with a single amplifier with variable gain was proposed in \cite{Ref17}.

On the other hand, the use of active components for amplification purposes inevitably increases the power consumption of the RIS, both due to the radio frequency power for signal amplification, and due to the additional hardware that is employed. Therefore, while it is clear that active RISs can increase the system spectral efficiency, it is not equally clear whether the same is true for the EE. Related works on the design of wireless networks aided by active RISs focus only on the system rate, e.g., \cite{Pan2022,PanActive,Liu2022}, neglecting the EE. This work aims at filling this gap, considering the maximization of the EE in a multi-user RIS-aided wireless network and comparing active and passive RISs. Moreover, the comparison is performed considering a RIS capable of global reflection, i.e. a recently-proposed kind of RIS in which the constraint on the reflected power is not applied to each reflecting element individually, but rather to the complete surface \cite{MDR22}. Finally, EE optimization is tackled  also with respect to the linear receive filters, which complicates the analysis.  
 
\section{System Model and Problem Formulation}\label{Sec:SystemModel}
\begin{figure}[!h]
\centering
\includegraphics[width=0.3\textwidth]{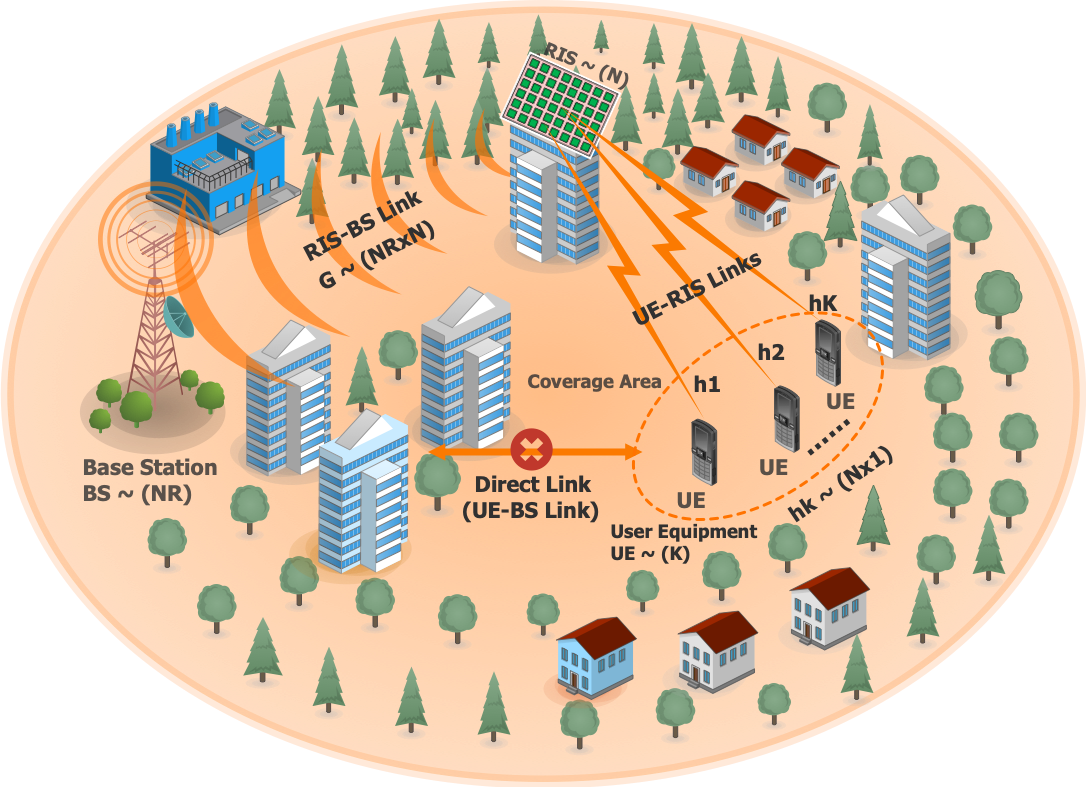}\caption{Considered wireless network.} \label{Fig:scenario}
\end{figure}

Let us consider the uplink of a multi-user system in which $K$ single-antenna mobile terminals communicate with a base station with $N_{R}$ antennas, through an RIS with $N$ reflecting elements (Fig. \ref{Fig:scenario}). 
Let us denote by $\bh_{k}$ the $N\times 1$ channel from user $k$ to the RIS, by $\bG$ the $N_{R}\times N$ channel from the RIS to the base station, by $\bGamma=\text{diag}(\gamma_{1},\ldots,\gamma_{N})$, the matrix whose diagonal contains the $N$ RIS reflection coefficients, and by $p_{k}$ and $s_{k}$ the transmit power and information symbol of user $k$. 

As already mentioned, most related works in the literature consider passive RISs with local reflection capabilities, i.e. they assume that the modulus of each RIS element is individually constrained to be less than a threshold, namely $|\gamma_{n}|^{2}\leq P_{R}$ for all $n=1,\ldots,N$, with $P_{R}\leq 1$. Instead, a more general scenario is that of active RISs with global reflection constraints, i.e. a constraint on the total power reflected by all of the RIS elements, and each element of the RIS is capable of amplifying the incoming signal thanks to the use of reflection-type amplifier implemented by different technologies, e.g. current inverters, current mirrors, or integrated circuits. 
In this case, the total power $P_{out}$ departing from the RIS can be larger than the incoming power $P_{in}$, with the difference being the radio-frequency power consumed by the RIS for amplification purposes. Then, the power constraint for an active RIS is written as $P_{out}-P_{in}\leq P_{R,max}$, with $P_{R,max}$ the maximum power that can be provided by the RIS for signal amplification. 
Mathematically speaking, the signals that arrive at and depart from the RIS can be written as $\by_{in}^{RIS}\!=\!\sum_{k=1}^{K}\sqrt{p}_{k}\bh_{k}s_{k}+\bn$ and $\by_{out}^{RIS}\!=\!\sum_{k=1}^{K}\sqrt{p}_{k}\bGamma\bh_{k}s_{k}+\bGamma\bn\!=\!\sum_{k=1}^{K}\sqrt{p}_{k}\bH_{k}\bgamma s_{k}+\bGamma\bn$, where $\bGamma\bh_{k}=\bH_{k}\bgamma$ with $\bH_{k}$ a diagonal matrix having the vector $\bh_{k}$ as its diagonal, i.e., $\bH_{k}=\textrm{diag}\left(h_{k}(1),\ldots,h_{k}(N)\right)$, $\bgamma=[\bGamma(1,1),\ldots,\bGamma(N,N)]^{T}$ contains the RIS coefficients, $p_{k}$ and $s_{k}$ are the $k$-th user's transmit power and information symbol, 
and $\bn$ is the thermal noise introduced by the transmit amplifier of the RIS\footnote{Note that due to the presence of power amplifiers, an active RIS introduces a noise amplification effect.}, modeled as a zero-mean, circularly symmetric Gaussian random variable with variance $\sigma_{\text{RIS}}^{2}$. The signal received at the final destination is $\br=\sum_{k=1}^{K}\sqrt{p}_{k}s_{k}\bG\bH_{k}\bgamma+\bG\bGamma\bn+\bw$, which, after applying the linear filter $\bc_{k}$, yields the SINR for user $k$
\begin{align}\label{Eq:SINR_a}
\text{SINR}_{k}&=\frac{p_{k}|\bc_{k}^{H}\bA_{k}\bgamma|^{2}}{\bc_{k}^{H}\bW\bc_{k}+\sum_{m\neq k}p_{m}|\bc_{k}^{H}\bA_{m}\bgamma|^{2}}\;,
\end{align}
with $\bA_{k}=\bG\bH_{k}$ for all $k$ and $\bW=\sigma^{2}\bI_{N_{R}}+\sigma_{\text{RIS}}^{2}\bG\bGamma\bGamma^{H}\bG^{H}$ is the covariance matrix of the overall colored noise at the receiver. Then, the power $P_{in}$ that arrives at the RIS and the power $P_{out}$ that departs from the RIS are given by
\begin{align}
P_{in}&\!\!=\!\!\sum_{k=1}^{K}\!p_{k}\|\bh_{k}\|^{2}\!\!+\!\sigma_{\text{RIS}}^{2}N\!=\!\!\!\sum_{k=1}^{K}p_{k}\tr(\bH_{k}\bH_{k}^{H})\!+\!\sigma_{\text{RIS}}^{2}N\\
P_{out}&=\sum_{k=1}^{K}p_{k}\tr(\bH_{k}\bgamma \bgamma^{H}\bH_{k}^{H})+\sigma_{\text{RIS}}^{2}\tr(\bgamma \bgamma^{H})
\end{align}
Then, the radio-frequency power consumed by the RIS is
\begin{align}\label{Eq:RF_Power_Active}
P_{out}\!-\!P_{in}&\!=\!\!\!\sum_{k=1}^{K}\!p_{k}\tr(\bH_{k}(\bgamma\bgamma^{H}\!\!\!-\!\bI_{N})\bH_{k}^{H})\!+\!\sigma_{\text{RIS}}^{2}\tr(\bgamma\bgamma^{H}\!\!\!-\!\bI_{N})\notag\\
&\hspace{-1cm}=\sum_{k=1}^{K}\!\left[p_{k}\tr(\bH_{k}(\bgamma\bgamma^{H}\!\!\!-\!\bI_{N})\bH_{k}^{H})\!+\!\frac{\sigma_{\text{RIS}}^{2}}{K}\tr(\bgamma\bgamma^{H}\!\!\!-\!\bI_{N})\!\right]\notag\\
&\hspace{-1cm}=\tr\left(\left(\bgamma\bgamma^{H}-\bI_{N}\right)\bR\right)\;,
\end{align}
wherein we have defined the positive definite, diagonal matrix 
\beq\label{Eq:R}
\bR=\sum_{k=1}^{K}\!\left(p_{k}\bH_{k}^{H}\bH_{k}\!+\!\frac{\sigma_{\text{RIS}}^{2}}{K}\bI_{N}\!\right)\!=\!\sum_{k=1}^{K}\!p_{k}\bH_{k}^{H}\bH_{k}\!+\!\sigma_{\text{RIS}}^{2}\bI_{N}.
\eeq
Thus, when the RIS operates in the active regime, it must hold $P_{out}-P_{in}\geq 0$, which yields the condition $\tr\left(\left(\bgamma\bgamma^{H}-\bI_{N}\right)\bR\right)\geq 0$, which, in turn, can be readily expressed as $\bgamma^{H}\bR\bgamma=\tr\left(\bR\bgamma\bgamma^{H}\right)\geq\tr\left(\bR\right)$. Next, the total power consumption of the system is given by 
\begin{align}
P_{tot}&=\textstyle\sum_{k=1}^{K}\mu_{k}p_{k}+p_{k}\tr\left(\bH_{k}(\bgamma\bgamma^{H}-\bI_{N})\bH_{k}^{H}\right)\notag\\
&+\sigma_{\text{RIS}}^{2}\tr\left(\bgamma\bgamma^{H}-\bI_{N}\right)+P_{c}\\
&=\textstyle\tr(\bR\bgamma\bgamma^{H})-\sigma_{\text{RIS}}^{2}N+\sum_{k=1}^{K}p_{k}(\mu_{k}-\|\bh_{k}\|^{2})+P_{c}\notag\;,
\end{align}
with $\mu_{k}$ the inverse of the efficiency of user's $k$ transmit amplifier, while $P_{c}=P_{0}+NP_{c,n}+P_{0,RIS}$, with $P_{c,n}$ the static power consumed by each active RIS element to enable the reconfiguration, and $P_{0,RIS}$ accounting for all other sources of static power consumption of the active RIS. Thus, the resulting global EE in the active RIS scenario is given by
\beq
\text{GEE}\!=\!\frac{B\sum_{k=1}^{K}\log_{2}\left(1+\text{SINR}_{k}\right)}{\tr(\bR\bgamma\bgamma^{H})\!-\!\sigma_{\text{RIS}}^{2}N\!+\!\sum_{k=1}^{K}p_{k}(\mu_{k}\!-\!\|\bh_{k}\|^{2})\!+\!P_{c}},
\eeq
with $B$ the communication bandwidth. Finally, the GEE maximization problem in the case of an active RIS with global reflection constraints is written as
\begin{subequations}\label{Prob:ActiveRIS}
\begin{align}
&\ds\max_{\bgamma,\bp,\bC}\; \text{GEE}(\bgamma,\bp,\bC)\label{Prob:aActiveRIS}\\
&\;\text{s.t.}\;\tr\left(\bR\right)\leq \bgamma^{H}\bR\bgamma\leq P_{R,max}+\tr\left(\bR\right)\label{Prob:bActiveRIS}\\
&\;\quad\;\;0\leq p_{k}\leq P_{max,k}\;\forall\;k=1,\ldots,K\;,\label{Prob:dActiveRIS}
\end{align}
\end{subequations}
wherein \eqref{Prob:bActiveRIS} ensures that the RIS operates in the active regime and that the RIS power budget is not exceeded. Also, observe that \eqref{Prob:bActiveRIS} is always feasible\footnote{Indeed, since $\bR$ is a diagonal matrix, the left-hand-side of \eqref{Prob:bActiveRIS} can be always fulfilled by choosing $|\gamma_{n}|=1$ for all $n=1,\ldots,N$, which falls back into the passive RIS regime.}. In order to tackle the non-convex Problem \eqref{Prob:ActiveRIS}, two methods are presented in Sections \ref{Sec:Design1} and \ref{Sec:Design2}. In both cases, it is assumed that channel realizations are reliably estimated and available for resource allocation. 

\section{First proposed approach}\label{Sec:Design1}
The first optimization method is based on the alternating optimization algorithm applied to the 
variables $\bgamma$, $\bp$, and $\bC$. 
\subsubsection{Optimization of $\bgamma$}
With respect to $\bgamma$, the problem is
\begin{subequations}\label{Prob:MaxGEEMF_gamma_A}
\begin{align}
&\ds\max_{\bgamma}\,\frac{\sum_{k=1}^{K}\log_{2}\left(1+\frac{p_{k}|\bc_{k}^{H}\bA_{k}\bgamma|^{2}}{\bc_{k}^{H}\bW\bc_{k}+\sum_{m\neq k}p_{m}|\bc_{k}^{H}\bA_{m}\bgamma|^{2}}\right)}{\tr(\bR\bgamma\bgamma^{H})+P_{c,eq}}\label{Prob:aMaxGEEMF_gamma_A}\\
&\;\text{s.t.}\;\tr{(\bR)}\leq\bgamma^{H}\bR\bgamma\leq P_{R,max}+\tr(\bR)
\label{Prob:bMaxGEEMF_gamma_A}
\end{align}
\end{subequations}
where $P_{c,eq}=\sum_{k}p_{k}(\mu_{k}-\|\bh_{k}\|^{2})+P_{c}-\sigma_{\text{RIS}}^{2}N$, and we recall that $\bW$ depends on the RIS vector $\bgamma$. Problem \eqref{Prob:MaxGEEMF_gamma_A} is more challenging than EE optimization with passive RISs, because: 1) the RIS vector $\bgamma$ appears also at the denominator of the objective \eqref{Prob:aMaxGEEMF_gamma_A}, which makes Problem \eqref{Prob:MaxGEEMF_gamma_A} a fractional program; 2) the RIS vector $\bgamma$ appears also in the receiver noise covariance matrix $\bW$; 3) the first inequality in \eqref{Prob:bMaxGEEMF_gamma_A} is a non-convex constraint. 

In order to deal with the fractional nature of Problem \eqref{Prob:MaxGEEMF_gamma_A}, the framework of fractional programming can be employed \cite{ZapNow15}. However, it can not be used directly, because the numerator of \eqref{Prob:aMaxGEEMF_gamma_A} is not a concave function. This issue will be addressed by resorting to the framework of sequential programming \cite{SeqCvxProg78}, and in particular to the sequential fractional programming method \cite{ZapNow15}. To begin with, let us express the numerator of \eqref{Prob:aMaxGEEMF_gamma_A}, and in particular the term $\bc_{k}^{H}\bW\bc_{k}$, as a function of the vector $\bgamma$ rather than the matrix $\bGamma$. To this end, defining $\bu_{k}=\bG^{H}\bc_{k}$, $\widetilde{\bU}_{k}=\text{diag}(|u_{1}|^{2},\ldots,|u_{N}|^{2})$, and plugging the expression of $\bW$ we obtain $\bc_{k}^{H}\bW\bc_{k}=\sigma^{2}\|\bc_{k}\|^{2}+\sigma_{\text{RIS}}^{2}\bgamma^{H}\widetilde{\bU}_{k}\bgamma$. Next, in order to apply the sequential fractional programming method, a concave lower-bound of the numerator of \eqref{Prob:aMaxGEEMF_gamma_A} must be found. To this end we apply the bound $\log_{2}\left(1\!+\!\frac{x}{y}\right)\!\!\geq\!\! \log_{2}\left(1\!+\!\frac{\bar{x}}{\bar{y}}\right)\!+\!\frac{\bar{x}}{\bar{y}}\left(\frac{2\sqrt{x}}{\sqrt{\bar{x}}}\!-\!\frac{x+y}{\bar{x}\!+\!\bar{y}}\!-\!1\right)$, which holds for any $x$, $y$, $\bar{x}$ and $\bar{y}$, and holds with equality whenever $x=\bar{x}$ and $y=\bar{y}$. Indeed, denoting by $\bar{\bgamma}$ any feasible vector of RIS reflection coefficients and applying the bound to each summand of the numerator of \eqref{Prob:aMaxGEEMF_gamma_A}, we obtain
\begin{align}
&R_{k}\!=\!\log_{2}\left(1+\text{SINR}_{k}\right)\label{Eq:Rbar1}\\
&\geq \!\log_{2}\left(\!1\!+\!\frac{p_{k}|\bc_{k}^{H}\bA_{k}\bar{\bgamma}|^{2}}{\sigma^{2}\|\bc_{k}\|^{2}\!+\!\sigma_{\text{RIS}}^{2}\bar{\bgamma}^{H}\widetilde{\bU}_{k}\bar{\bgamma}\!+\!\sum_{m\neq k}p_{m}|\bc_{k}^{H}\bA_{m}\bar{\bgamma}|^{2}}\right)\notag\\
&\!+\!\frac{p_{k}|\bc_{k}\bA_{k}\bar{\bgamma}|^{2}}{\sigma^{2}\|\bc_{k}\|^{2}\!+\!\sigma_{\text{RIS}}^{2}\bar{\bgamma}^{H}\widetilde{\bU}_{k}\bar{\bgamma}\!+\!\sum_{m\neq k}p_{m}|\bc_{k}^{H}\bA_{m}\bar{\bgamma}|^{2}}\!\left(\!\frac{2|\bc_{k}^{H}\bA_{k}\bgamma|}{|\bc_{k}^{H}\bA_{k}\bar{\bgamma}|}\right.\notag\\
&\left.-\frac{\sigma^{2}\|\bc_{k}\|^{2}\!+\!\sigma_{\text{RIS}}^{2}\bgamma^{H}\widetilde{\bU}_{k}\bgamma\!+\!\sum_{m=1}^{K}p_{m}|\bc_{k}^{H}\bA_{m}\bgamma|^{2}}{\sigma^{2}\|\bc_{k}\|^{2}\!+\!\sigma_{\text{RIS}}^{2}\bar{\bgamma}^{H}\widetilde{\bU}_{k}\bar{\bgamma}\!+\!\sum_{m=1}^{K}p_{m}|\bc_{k}^{H}\bA_{m}\bar{\bgamma}|^{2}}\!-\!1\right)\!=\!\bar{R}_{k}\notag
\end{align}
Defining $\bar{L}_{k}\!=\!\sigma^{2}\|\bc_{k}\|^{2}\!+\!\sigma_{\text{RIS}}^{2}\bar{\bgamma}^{H}\widetilde{\bU}_{k}\bar{\bgamma}\!+\!\sum_{m\neq k}p_{m}|\bc_{k}^{H}\bA_{m}\bar{\bgamma}|^{2}$,  $\bar{A}_{k}\!=\!\log_{2}\left(1+p_{k}|\bc_{k}^{H}\bA_{k}\bar{\bgamma}|^{2}/\bar{L}_{k}\right)$, $\bar{B}_{k}=p_{k}|\bc_{k}^{H}\bA_{k}\bar{\bgamma}|^{2}/\bar{L}_{k}\;,\;\bar{D}_{k}=2/|\bc_{k}^{H}\bA_{k}\bar{\bgamma}|$, $\bar{E}_{k}=1/\left(\bar{L}_{k}+p_{k}|\bc_{k}^{H}\bA_{k}\bar{\bgamma}|^{2}\right)$, and $\bar{F}_{k}=\bar{E}_{k}\sigma^{2}\|\bc_{k}\|^{2}+1$, 
we can write $\bar{R}_{k}$ as
\begin{align}\label{Eq:BarR_Active}
\bar{R}_{k}&=\bar{A}_{k}+\bar{B}_{k}\Bigg(\bar{D}_{k}|\bc_{k}^{H}\bA_{k}\bgamma|-\bar{F}_{k}\notag\\
&-\bar{E}_{k}\left(\sigma_{\text{RIS}}^{2}\bgamma^{H}\widetilde{\bU}_{k}\bgamma+\textstyle\sum_{m=1}^{K}p_{m}|\bc_{k}^{H}\bA_{m}\bgamma|^{2}\right)\Bigg)
\end{align}
Let us observe that the term $-\sigma_{\text{RIS}}^{2}\bgamma^{H}\widetilde{\bU}_{k}\bgamma$ is concave and thus the only non-concave term in $\bar{R}_{k}$ is $|\bc_{k}^{H}\bA_{k}\bgamma|$. However, since $|\bc_{k}^{H}\bA_{k}\bgamma|$ is convex in $\bgamma$, it is lower-bounded by its first-order Taylor expansion around any point $\bar{\bgamma}$, which yields the following  concave lower-bound of the generic summand in the numerator of \eqref{Prob:aMaxGEEMF_gamma_A}  
\begin{align}
R_{k}&\geq \bar{R}_{k}\geq \bar{A}_{k}+\bar{B}_{k}\Bigg(\bar{D}_{k}\Bigg(|\bc_{k}^{H}\bA_{k}\bar{\bgamma}|\\
&+\Re\left\{\frac{\bA_{k}^{H}\bc_{k}\bc_{k}^{H}\bA_{k}\bar{\gamma}}{|\bc_{k}^{H}\bA_{k}\bar{\bgamma}|}(\bgamma-\bar{\bgamma})\right\}\Bigg)\\
&-\bar{E}_{k}\left(\sigma_{\text{RIS}}^{2}\bgamma^{H}\widetilde{\bU}_{k}\bgamma+\!\!\!\sum_{m=1}^{K}p_{m}|\bc_{k}^{H}\bA_{m}\bgamma|^{2}\right)-\bar{F}_{k}\Bigg)\!=\!\widetilde{R}_{k}\notag
\end{align}
Finally, it remains to deal with the first, non-convex inequality constraint in \eqref{Prob:bMaxGEEMF_gamma_A}. This is accomplished by employing the sequential approximation framework. Specifically, since $\bgamma^{H}\bR\bgamma$ is a convex function of $\bgamma$, it is lower-bounded by its first-order Taylor expansion around any point $\bar{\gamma}$, i.e., $\bgamma^{H}\bR\bgamma\geq \bar{\bgamma}^{H}\bR\bar{\bgamma}+2\Re\{\bar{\bgamma}^{H}\bR(\bgamma-\bar{\bgamma})\}$. Thus, in each iteration of the sequential method, the problem to be solved is 
\begin{subequations}\label{Prob:GEE_Gamma_Active}
\begin{align}
&\ds\max_{\bgamma}\frac{\sum_{k=1}^{K}\widetilde{R}_{k}(\bgamma)}{\tr(\bR\bgamma\bgamma^{H})+P_{c,eq}}\label{Prob:aGEE_Gamma_Active}\\
&\;\text{s.t.}\;\bgamma^{H}\bR\bgamma\leq P_{R,max}+\tr{(\bR)}\;.\label{Prob:bGEE_Gamma_Active}\\
&\;\quad\;\bar{\bgamma}^{H}\bR\bar{\bgamma}+2\Re\{\bar{\bgamma}^{H}\bR(\bgamma-\bar{\bgamma})\}\geq \tr(\bR)\;.\label{Prob:cGEE_Gamma_Active}
\end{align}
\end{subequations}
The objective \eqref{Prob:aGEE_Gamma_Active} has a concave numerator and a convex denominator, since $\widetilde{R}_{k}$ is concave and $\bR$ is positive definite, while \eqref{Prob:cGEE_Gamma_Active} is a linear constraint. Thus, \eqref{Prob:GEE_Gamma_Active} is a pseudo-concave maximization subject to convex constraints, which can be solved by fractional programming methods. 
\subsubsection{Optimization of $\bp$}
Defining $a_{k,m}=|\bc_{k}^{H}\bA_{m}\bgamma|^{2}$, for all $m$ and $k$, $d_{k}=\bc_{k}^{H}\bW\bc_{k}$, $P_{c,eq}=\tr(\bR\bgamma\bgamma^{H})\!-\!\sigma_{\text{RIS}}^{2}N+P_{c}$, and $\mu_{k,eq}=\mu_{k}\!-\!\|\bh_{k}\|^{2}$, the power optimization problem is
\begin{subequations}\label{Prob:MaxGEEMF_power}
\begin{align}
&\ds\max_{\bp}\,\frac{\sum_{k=1}^{K}\ds\log_{2}\left(1+\frac{p_{k}a_{k,k}}{d_{k}+\sum_{m\neq k}p_{m}a_{k,m}}\right)}{\sum_{k=1}^{K}\mu_{k,eq}p_{k}+P_{c,eq}}\label{Prob:aMaxGEEMF_power}\\
&\;\text{s.t.}\;0\leq p_{k}\leq P_{max,k}\;,\forall\; k=1,\ldots,K
\end{align}
\end{subequations}
Since the numerator of \eqref{Prob:aMaxGEEMF_power} is not a concave function of $\bp$, the objective in \eqref{Prob:aMaxGEEMF_power} is not a pseudo-concave function and thus it is computationally unfeasible to solve \eqref{Prob:MaxGEEMF_power} by fractional programming  \cite{ZapNow15}. Then, we resort to the sequential fractional programming method \cite{ZapNow15}, in order to derive a pseudo-concave lower-bound of \eqref{Prob:aMaxGEEMF_power}, which can be maximized by fractional programming. To this end, let us express \eqref{Prob:aMaxGEEMF_power} as  
\begin{align}
\text{GEE}(\bp)&=\underbrace{\frac{\sum_{k=1}^{K}\log_{2}\left(\!d_{k}+\sum_{k=1}^{K}p_{k}a_{k,k}\right)}{\sum_{k=1}^{K}\mu_{k,eq}p_{k}+P_{c,eq}}}_{\ds g_{1}(\bp)}\\
&-\underbrace{\frac{\sum_{k=1}^{K}\log_{2}\left(d_{k}+\sum_{m\neq k}p_{m}a_{k,m}\right)}{\sum_{k=1}^{K}\mu_{k,eq}p_{k}+P_{c,eq}}}_{\ds g_{2}(\bp)}.\notag
\end{align}
Since the numerator of $g_{2}(\bp)$ is concave, a pseudo-concave lower-bound of $\text{GEE}(\bp)$, denoted by $\widetilde{\text{GEE}}(\bp)$, is obtained replacing the numerator of $g_{2}(\bp)$ with its first-order Taylor expansion around any feasible point $\bar{p}$ \cite{ZapNow15}.
In each iteration of the sequential method, the function $\widetilde{\text{GEE}}(\bp)$ is maximized, with the constraint $p_{k}\in[0,P_{max,k}]$, for all $k=1,\ldots,K$, which is efficiently solved by fractional programming.
\subsubsection{Optimization of $\bC$}\label{Sec:C}
The optimization of the receive filters in $\bC$ affects only the numerator of the GEE. Moreover, it can be decoupled over the users, thus reducing to maximizing the individual rate of the users. The solution to this problem is well-known to be the linear MMSE receiver, which for the case at hand, is expressed according to the formula $\bc_{k}\!=\!\sqrt{p}_{k}\bM_{k}^{-1}\bA_{k}\bgamma$, with $\bM_{k}=\sum_{m\neq k}p_{m}\bA_{m}\bgamma\bgamma^{H}\bA_{m}^{H}+\bW$. 

\section{Second proposed approach}\label{Sec:Design2}
While the previous approach considers the alternating optimization of three variables, namely, the transmit powers $\bp$, the receive filters $\bc_{1},\ldots,\bc_{K}$, and the RIS reflection coefficients $\bgamma$, a different strategy is that of embedding the optimal expression of the linear MMSE filters into the GEE, and optimizing the resulting expression with respect to $\bp$ and $\bgamma$. While this is expected to yield better performance, an intuition that will be confirmed by numerical results, it also leads to a more challenging optimization problem, due to the more involved expression of the EE function. To elaborate, accounting for the fact that RIS amplification colors the noise at the receiver, the sum-rate with linear MMSE filtering is written as 
\begin{align}
\text{SR}_{\text{MMSE}}&=\textstyle\sum_{k=1}^{K}\log_{2}\left|\bW+\textstyle\sum_{m=1}^{K}p_{m}\bA_{m}\bgamma\bgamma^{H}\bA_{m}^{H}\right|\notag\\
&-\textstyle\sum_{k=1}^{K}\log_{2}\left|\bW+\textstyle\sum_{m\neq k}p_{m}\bA_{m}\bgamma\bgamma^{H}\bA_{m}^{H}\right|\label{Eq:SR_MMSE}
\end{align}
Thus, the GEE maximization problem can be formulated as
\begin{subequations}\label{Prob:GEEMaxActive}
\begin{align}
&\ds\max_{\bgamma,\bp}\;\frac{\text{SR}_{\text{MMSE}}(\bgamma,\bp)}{\tr(\bR\bgamma\bgamma^{H})-\sigma_{\text{RIS}}^{2}N+\sum_{k=1}^{K}p_{k}(\mu_{k}-\|\bh_{k}\|^{2})+P_{c}}\label{Prob:aGEEMax}\\
&\;\text{s.t.}\;\tr(\bR)\leq \tr(\bR\bgamma\bgamma^{H})\leq P_{R,max}+\tr(\bR)\label{Prob:bGEEMax}\\
&\;\quad\;p_{k}\in[0,P_{max,k}]\;,\forall\;k=1,\ldots,K\label{Prob:dGEEMax}\;,
\end{align}
\end{subequations}
where we have exploited that $\bgamma^{H}\bR\bgamma=\tr(\bR\bgamma\bgamma^{H})$. Problem \eqref{Prob:GEEMaxActive} can be tackled by alternating optimization of $\bgamma$ and $\bp$ as discussed in the rest of this section. 

\subsubsection{Optimization of $\bgamma$}
The problem can be stated as 
\begin{subequations}\label{Prob:GEEMaxActive_2}
\begin{align}
&\ds\max_{\bgamma}\;\frac{\text{SR}_{\text{MMSE}}(\bgamma)}{\tr(\bR\bgamma\bgamma^{H})+P_{c,eq}}\label{Prob:aGEEMaxActive_2}\\
&\;\text{s.t.}\;\tr(\bR)\leq\tr(\bR\bgamma\bgamma^{H})\leq P_{R,max}+\tr(\bR)\label{Prob:bGEEMaxActive_2}
\end{align}
\end{subequations}
To begin with, let us observe that the noise covariance matrix $\bW$ can be expressed as $
\bW=\sigma^{2}\bI_{N_{R}}+\sigma_{\text{RIS}}^{2}\bG\bGamma\bGamma^{H}\bG^{H}=\sigma^{2}\bI_{N_{R}}+\sigma_{\text{RIS}}^{2}\bG\text{diag}(\bgamma\bgamma^{H})\bG^{H}$.
Then, defining $\bX=\bgamma\bgamma^{H}$, the numerator of the objective in \eqref{Prob:aGEEMaxActive_2} can be expressed as
\begin{align}\label{Eq:SR_X_Active}
&\text{SR}_{\text{MMSE}}(\bX)=\\
&\underbrace{\sum_{k=1}^{K}\!\log_{2}\!\left|\sigma^{2}\bI_{N_{R}}\!\!+\!\sigma_{\text{RIS}}^{2}\bG\text{diag}(\bX)\bG^{H}\!\!+\!\sum_{m=1}^{K}p_{m}\bA_{m}\bX\bA_{m}^{H}\right|}_{\ds G_{1}(\bX)}\!-\notag\\
&\underbrace{\sum_{k=1}^{K}\!\log_{2}\left|\sigma^{2}\bI_{N_{R}}\!\!+\!\sigma_{\text{RIS}}^{2}\bG\text{diag}(\bX)\bG^{H}\!\!+\!\sum_{m\neq k}p_{m}\bA_{m}\bX\bA_{m}^{H}\right|}_{\ds G_{2}(\bX)}\notag,
\end{align}
where $\text{diag}(\bX)$ is the diagonal of the matrix $\bX$, which is a linear function of $\bX$.\footnote{Observe also that $\text{diag}(\bX)=\bX\odot\bI_{N}$, with $\odot$ denoting the component-wise product.} Thus, being the difference of two concave functions, a concave lower-bound of \eqref{Eq:SR_X_Active} is found linearizing the negative sum in \eqref{Eq:SR_X_Active} around any point $\bar{\bX}$. 

\noindent Specifically, it holds 
\begin{align}
\text{SR}_{\text{MMSE}}(\bX)&\geq G_{1}(\bX)-G_{2}(\bar{\bX})\\
&-\Re\!\left\{\tr\!\left(\nabla G_{2}(\bar{\bX})^{H}(\bX-\bar{\bX})\right)\right\}=\widetilde{\text{SR}}_{\text{MMSE}}(\bX)\;,\notag
\end{align}
with 
\beq
\nabla G_{2}(\bX)\!\!=\!\!\!\sum_{k=1}^{K}\!\!\left(\!\sigma_{\text{RIS}}^{2}(\bG^{H}\bT_{k}^{-1}\bG)\!\odot\!\bI_{N}\!\!+\!\!\!\!\sum_{m\neq k}\!\!p_{m}\bA_{m}^{H}\bT_{k}^{-1}\bA_{m}\!\!\right)\!\!,
\eeq
and $\bT_{k}\!=\!\sigma^{2}\bI_{N_{R}}\!+\!\sigma_{\text{RIS}}^{2}\bG\text{diag}(\bX)\bG^{H}+\sum_{m\neq k}p_{m}\bA_{m}\bX\bA_{m}^{H}$.
Then, the problem to be solved in each iteration of the sequential method can be formulated as 
\begin{subequations}\label{Prob:MaxGEEMMMSE_Xrel_SQ}
\begin{align}
&\ds\max_{\bX\succeq \bzero}\,\frac{\widetilde{\text{SR}}_{\text{MMSE}}(\bX)}{\tr(\bR\bX)+P_{c,eq}}\label{Prob:aMaxGEEMMSE_Xrel_SQ}\\
&\;\text{s.t.}\;\tr(\bR)\leq \tr(\bR\bX)\leq P_{R,max}+\tr(\bR)\label{Prob:bMaxGEEMMSE_Xrel_SQ}\;,
\end{align}
\end{subequations}
where we have lifted the rank-one constraint on $\bX$, thus employing the semidefinite relaxation method \cite{LuoSDR}. It can be seen that \eqref{Prob:aMaxGEEMMSE_Xrel_SQ} has a concave numerator and a convex denominator, and thus is a pseudo-concave function. Then, since the constraints of Problem \eqref{Prob:MaxGEEMMMSE_Xrel_SQ} are affine, the problem can be globally solved by standard fractional programming tools. Upon convergence of the sequential procedure, if the convergence point $\bX^{*}$ has unit-rank, then $\bX^{*}$ is also feasible for the original problem. Otherwise, a feasible solution can be obtained by randomization techniques \cite{LuoSDR}. 

\subsubsection{Optimization of $\bp$} 
The problem can be stated as
\begin{subequations}\label{Prob:MaxGEEMMMSE_power}
\begin{align}
&\ds\max_{\bp}\,\text{GEE}_{\text{MMSE}}=\frac{\text{SR}_{\text{MMSE}}(\bp)}{\sum_{k=1}^{K}\mu_{k,eq}p_{k}+P_{c,eq}}\label{Prob:aMaxGEEMMSE_power}\\
&\;\text{s.t.}\;p_{k}\in[0, P_{max,k}]\;\forall\;k=1,\ldots,K\;,
\end{align}
\end{subequations}
which can be tackled by sequential fractional programming by linearizing the negative term in the numerator of \eqref{Prob:aMaxGEEMMSE_power}. Then, in each iteration of the sequential fractional programming algorithm, the surrogate problem to be solved is
\begin{align}\label{Prob:GEEMax_SFP_Active}
&\ds\max_{\bp}\widetilde{\text{GEE}}_{\text{MMSE}}(\bp),\text{s.t.}\,p_{k}\!\in\![0,\!P_{max,k}],\forall k\!=\!1,\ldots,K,
\end{align}
with the objective given by 
\begin{align}
&\widetilde{\text{GEE}}_{\text{MMSE}}(\bp)\!=\!\frac{\ds\sum_{k=1}^{K}\log_{2}\left|\bW\!+\!\textstyle\sum_{m=1}^{K}p_{m}\bA_{m}\bgamma\bgamma^{H}\bA_{m}^{H}\right|}{\ds\sum_{k=1}^{K}\mu_{k,eq}p_{k}+P_{c,eq}}\notag\\
&-\frac{\ds\sum_{k=1}^{K}\log_{2}\left|\bW\!+\!\!\!\sum_{m\neq k}\bar{p}_{m}\bA_{m}\bgamma\bgamma^{H}\bA_{m}^{H}\right|}{\ds\sum_{k=1}^{K}\mu_{k,eq}p_{k}+P_{c,eq}}\!-\!\frac{\left(\nabla_{\bp} F(\bar{\bp})\right)^{T}\!\!\!\left(\bp\!-\!\bar{\bp}\right)}{\ds\sum_{k=1}^{K}\mu_{k,eq}p_{k}+P_{c,eq}}\notag,
\end{align}
wherein $\bar{\bp}$ is the point around which the first-order Taylor expansion is computed, and, for any $i=1,\ldots,K$, it holds

\beq
\frac{\partial F}{\partial p_{i}}=\sum_{k\neq i}\tr\left(\left(\bW+\!\!\sum_{m\neq k}p_{m}\bA_{m}\bgamma\bgamma^{H}\bA_{m}^{H}\right)^{-1}\!\!\!\!\!\bA_{i}\bgamma\bgamma^{H}\bA_{i}^{H}\!\!\right)\;.\notag
\eeq
\subsubsection{Convergence and complexity}\label{Sec:Complexity}
Both methods are provably convergent, which follows from the convergence of the alternating optimization and sequential optimization methods \cite{ZapTSP18}. Moreover, they both enjoy polynomial complexity in the number of variables. Specifically, the computational complexity of the method from Section\footnote{We neglect the complexity of computing the receive matrix $\bC$, which can be accomplished by the closed-form expressions in Section \ref{Sec:C}.} \ref{Sec:Design1} is $\mathcal{C}=\mathcal{O}\left(I\left(I_{\gamma}(N+1)^{\alpha}+I_{p}(K+1)^{\beta}\right)\right)$, wherein $I_{\gamma}$ and $I_{p}$ are the number of iterations that the sequential method applied to $\gamma$ and $\bp$ requires to converge; $I$ is the number of iterations that the outer alternating optimization requires to converge; $(N+1)^{\alpha}$ is the complexity of each RIS optimization and $(K+1)^{\beta}$ is the complexity of each users' powers optimization\footnote{Pseudo-concave maximizations with $n$ variables can be restated as concave maximizations with $n+1$ variables \cite{ZapNow15}, and the complexity of a concave maximization is polynomial in the number of variables with exponent in the interval $[1,4]$. Thus, $\alpha$ and $\beta$ belong to the interval $[1,4]$ \cite{ComplexityBook}.}. A similar complexity analysis holds for the method from Section \ref{Sec:Design2}. Due to the involved mathematical structure of the problems at hand, the specific values of $\alpha$ and $\beta$, as well as of the iterations numbers $I_{\gamma}$, $I_{p}$, $I$, are not known in closed-form. Similarly, it is difficult to relate the choice of the initialization points to the convergence rate. Nevertheless, the theoretical convergence claim and complexity formulas remain valid for any choice of the initialization points. 

\section{Numerical Results}\label{Sec:Numerics}
We consider an instance of the multi-user network described in Section \ref{Sec:SystemModel}, with $K=4$, $N_{R}=4$, $N=100$, $B=20\,\textrm{MHz}$, $P_{0}=40\,\textrm{dBm}$, $P_{0,RIS}=30\,\textrm{dBm}$, $P_{c,n}=0\,\textrm{dBm}$, $P_{R,max}=10\,\textrm{dBW}$. The noise power spectral density is $-174\,\textrm{dBm/Hz}$, and a noise figure of $10\,\textrm{dB}$ has been considered. The mobile users are randomly placed in an area with radius $100\,\textrm{m}$ around the RIS, and the base station is placed $50\,\textrm{m}$ away from the RIS. The users have a random height in $[0,5]\,\textrm{m}$, while the RIS and base station have heights of $15\,\textrm{m}$ and $10\,\textrm{m}$, respectively. A power decay factor $\eta=4$ has been considered, while the fading component of all channels follows the Rice model, with factors $K_{t}=4$ for the channel from the RIS to the base station and $K_{r}=2$ for the channels from the mobile users to the RIS. 

\begin{figure}[!h]
\centering
\includegraphics[width=0.5\textwidth]{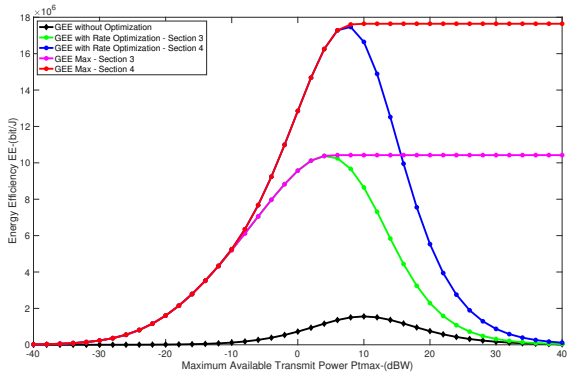}\caption{Achieved GEE versus $P_{max}$. $K=4$, $N_{R}=4$, $N=100$.} \label{fig:EEvsP}
\end{figure}

Fig. \ref{fig:EEvsP} shows the GEE achieved by:
\begin{itemize}
\item (a) the resource allocation obtained by maximizing the GEE through the method from Section \ref{Sec:Design1}.
\item (b) the resource allocation obtained by maximizing the GEE through the method from Section \ref{Sec:Design2}.  
\item (c) the resource allocation obtained by maximizing the spectral efficiency through the method from Section\footnote{Both the methods from Sections \ref{Sec:Design1} and \ref{Sec:Design2} can be be specialized to perform sum-rate maximization by simply setting $\mu_{k}=0$, for all $k=1,\ldots,K$.} \ref{Sec:Design1}. 
\item (d) the resource allocation obtained by maximizing the spectral efficiency through the method from Section \ref{Sec:Design2}. 
\item (e) uniform power allocation and random RIS phases . 
\end{itemize}
As anticipated, the method from Section \ref{Sec:Design2} significantly outperforms the approach from Section \ref{Sec:Design1}, thanks to the exploitation of the mathematical structure of the optimal receive filter, rather than simply updating it within the alternating maximization algorithm. Moreover, a huge gain is obtained compared to Case (e), in which the system radio resources are not optimized. As expected, the value of the GEE saturates after the point at which $P_{max}$ is large enough to allow reaching the maximum of the GEE function. Indeed, after this point, changing the resource allocation would only decrease the GEE. 

\begin{figure}[!h]
\centering
\includegraphics[width=0.5\textwidth]{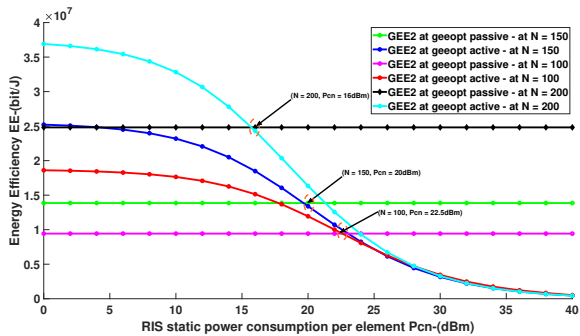}\caption{Achieved GEE versus $P_{c,n}$ for the active and passive RIS scenarios. $K=4$, $N_{R}=4$, $N=100$, $150$, $200$.} \label{fig:EEvsPcn}
\end{figure}

Fig. \ref{fig:EEvsPcn} compares the GEE obtained by the GEE-maximizing method from Section \ref{Sec:Design2} to that obtained by a passive RIS versus the power consumption $P_{c,n}$ of each active RIS element. As for the passive RIS scenario, the optimization is performed as a special case of the algorithm from Section \ref{Sec:Design2}, by setting $P_{R,max}=0$, which leads to setting $|\gamma_{n}|=1$ for all $n=1,\ldots,N$, and leaves the phases of the RIS reflection coefficients as optimization variables. Also, for the passive RIS case, fixed values of $P_{c,n}^{(p)}=0\,\textrm{dBm}$ and $P_{0,RIS}^{(p)}=20\,\textrm{dBm}$ have been considered. The results show the operating regime in which the active RIS is more energy-efficient than the passive one. Eventually, as  $P_{c,n}$ of each active RIS element increases, the passive RIS becomes more energy-efficient than the active one. Moreover, the crossing point is encountered for lower values of $P_{c,n}$ when the RIS is equipped with more active elements $N$. Thus, there is a trade-off between the GEE and the number of active RIS elements. 

Finally, Table \ref{Tab:1} analyzes the convergence time of the proposed methods. Let us denote by $T_{1,p}$ and $T_{2,p}$ the convergence time of the algorithms from Sections \ref{Sec:Design1} and \ref{Sec:Design2} when applied with a passive RIS. Similarly, let us denote by $T_{1,a}$ and $T_{2,a}$ the convergence time of the algorithms from Sections \ref{Sec:Design1} and \ref{Sec:Design2} when applied with an active RIS. Simulations have been performed by a workstation equipped with an AMD Ryzen 9, 5950x, 16-core processor. All convex optimization problems have been solved the CVX software package. Convergence is declared when the absolute error between the GEE values in two successive iterations is lower than $\epsilon=10^{-6}$. The results show that applying the proposed methods with active RIS can take longer than with a passive RIS. This is expected since the formulation of the active scenario leads to a more general and complex mathematical expression of the GEE and problem constraints. Similarly, the use of the method from Section \ref{Sec:Design2} can take longer than the one from Section \ref{Sec:Design1}, which is also expected due to the more involved expression of the objective function that we obtain by embedding the optimal expression of the receive filters into the SINR formula. This leads to a larger number of iterations required for convergence and, possibly, to larger values of the exponents $\alpha$ and $\beta$ in the computational complexity formulas derived in Section \ref{Sec:Complexity}. Thus, there is a trade-off between the first and second proposed methods. The former achieves a lower GEE value, but requires a lower computational complexity than the latter. 
\begin{table}[h]
\caption{Comparison between convergence time of the proposed algorithms with active and passive RIS.}\label{Tab:1}
\begin{center}\begin{tabular}{ccccc}
\hline
$P_{max}\,\textrm{[dBW]}$ & $T_{1,a}/T_{1,p}$ & $T_{2,a}/T_{2,p}$ & $T_{2,p}/T_{1,p}$ & $T_{2,a}/T_{1,a}$ \\
\hline
-40 & 2.26 & 1.88 & 12.58 & 10.45\\
-20 & 3.34 & 2.58 & 17.49 & 13.52\\
0 & 6.13 & 1,90 & 26.67 & 20.64\\
20 & 11.98 & 8.67 & 40.45 & 29.26\\
\hline\\
\end{tabular} 
\end{center}
\end{table} 

\section{Conclusions}\label{Sec:Conclusions}
This work has proposed two provably convergent algorithms with polynomial complexity for EE maximization in a wireless network aided by an active RIS. The RIS coefficients, the mobile users' transmit powers, and linear receive filters have been jointly optimized. Numerical results show the merits of the proposed algorithms and highlight a trade-off between the EE of active and passive RISs, in terms of the number of RIS elements and the additional power consumption due to the presence of the active-load hardware. 
\bibliographystyle{IEEEtran}
\bibliography{HuaweiBib,references,references_Active,FracProg}

\end{document}